\documentclass[a4paper,12pt]{article}

\newcommand{\di}{\displaystyle}
\newcommand{\ov}{\over}

\newcommand{\al}{\alpha}
\newcommand{\be}{\beta}

\newcommand{\na}{\nabla}

\newcommand{\ty}{\infty}

\newcommand{\pa}{\partial}
\newcommand{\va}{\varphi}
\newcommand{\ld}{\ldots}

\newcommand{\la}{\lambda}

\newcommand{\x}{{\rm \bf x}}
\newcommand{\g}{{\rm \bf g}}
\newcommand{\qed}{\hfill $\sqcap\!\!\!\!\sqcup$} 
\begin{document}
%
%
%
\title{Harmonicity and submanifold maps}
\author{ Constantin Udri\c ste, Vasile Arsinte and Andreea Bejenaru}
\date{}
\maketitle

\pagestyle{plain}
\begin{abstract}
{The aim of this paper is fourfold. Firstly, we introduce and study the $f$-ultra-harmonic maps.
Secondly, we recall the geometric dynamics generated by a first order normal PDE
system and we give original results regarding the geometric dynamics generated
by other first order PDE systems. Thirdly, we determine the Gauss PDEs and the
fundamental forms associated to integral manifolds of first order PDE systems.
Fourthly, we change the Gauss PDEs into a geometric dynamics on the jet bundle
of order one, showing that there exist an infinity of Riemannian metrics such that the lift of a
submanifold map into the first order jet bundle to be an ultra-potential map.

Keywords: harmonic map, ultra-potential map, generalized potential map, general harmonicity, Gauss equation

2010 Mathematics subject classification: 35J47, 35K70, 49J20}

\end{abstract}

\section{Generalized Harmonic and Potential Maps}

All maps throughout the paper are smooth, while manifolds are real,
finite-dimensional, Hausdorff, second-countable and connected.

Let $(N,h)$ be a Riemannian manifold of dimension $m$ and let $M$ be differential manifold with dimension $n$.
Hereafter we shall assume that the manifold $N$ is oriented. Greek (Latin)
letters will be used for indexing the components of geometrical objects
attached to the manifold $N$ (manifold $M$). Local coordinates on $N$ will be
written
$
t = (t^\al), \quad \al = 1, \ld, m,
$ and those on $M$ will be
$
x = (x^i), \quad i = 1, \ld, n.
$
The components of the corresponding metric tensor $h$ and Christoffel symbols on the manifold $N$
will be denoted by $h_{\al \beta}, H^\al_{\beta \gamma}$.

The product manifold $N \times M$ is endowed with the coordinates $(t^\al, x^i)$ and the
first order jet manifold $J^1 (N,M)$, called the configuration bundle, is endowed
with the adapted coordinates $(t^\al, x^i, x^i_\al)$. The distinguished tensors fields
and other distinguished geometrical objects on $N \times M$ are introduced using
the geometry of the jet bundle $J^1 (N,M)$ [4], [5], [15].

Let $\va:N \to M$, $\va (t) = x$, $x^i = x^i (t^\al)$ be a $C^\ty$
map (parameterized $m$-sheet). For a fixed symmetric (possible
degenerated) $(0,2)$-tensor field $f=(f_{ij})$ on $M$, we attach the
{\it $f$-energy density Lagrangian} defined by
$$
E_f (\va) (t) = {1 \ov 2} h^{\al\beta} (t) f_{ij} (x(t)) x^i_\al (t)
x^j_\beta (t)
\leqno(1)$$
and the total energy
$$
E_f (\va, \Omega) = \int_\Omega E_f (\va) (t) dv_h,
$$
where $|h|=\mbox{det}\,h$ and $dv_h = \sqrt{|h|}\, dt^1 \wedge \ld \wedge dt^m$ denotes the volume
element induced by the Riemannian metric $h$.

 {\bf Definition 1.} A map $\va$ is called {\it $f$ - ultra-harmonic map} if it is a critical point for the $f$-energy functional $E_f$,
i.e., an extremal of the Lagrangian
$$
L_1 = E_f (\va) (t) \sqrt{|h|},
$$
for all compactly supported variations.

If we denote by
$$F_{jk|i}=\frac{1}{2}\left\{\frac{\pa f_{ij}}{\pa x^k}+\frac{\pa f_{ik}}{\pa x^j}-\frac{\pa f_{jk}}{\pa x^i}\right\}$$
{\it the Christoffel symbols of the first type attached to tensor $f$} and if we introduce the distinguished tensor field
$$x_{i\al\beta} = f_{ij}{\pa^2 x^j \ov
\pa t^\al \pa t^\beta} - H^\gamma_{\al\beta}f_{ij} x^j_\gamma + F_{jk|i} x^j_\al x^k_\beta,$$
then an {\it $f$-ultra-harmonic map equation} is written in local coordinates as
$$h^{\al\beta}x_{i\alpha\beta} = 0 \leqno(2)$$
(a nonlinear {\it ultra-parabolic-hyperbolic PDE system} of second order).

Let $g$ be a Riemannian metric on the manifold $M$ and $G^k_{ij}$ be the corresponding Christoffel symbols.
In particular, if $f=g$, we obtain the definition of classical {\it harmonic maps} [1]-[7], [11]-[13], [15]-[17], [19]-[20]. Indeed,
the classical form of the {\it kinetic energy density} corresponding to the map $\va$ is
$$
E_g (\va) (t) = {1 \ov 2} h^{\al\beta} (t) g_{ij} (x(t)) x^i_\al (t)x^j_\beta (t)
$$
and the {\it harmonic map equation}  (a system of nonlinear {\it elliptic-Laplace PDEs} of second order),
is expressed in local coordinates by
$$
h^{\al\beta} x^i_{\alpha\beta} = 0,
$$
where $$x^i_{\al\beta} = {\pa^2 x^i \ov
\pa t^\al \pa t^\beta} - H^\gamma_{\al\beta} x^i_\gamma + G^i_{jk} x^j_\al x^k_\beta.$$

Let  $T = (T^i_\al (t))$ be a given $C^\ty$ distinguished tensor field on $N$, let $g$ be a
Riemannian structure on $M$ and $f$ be a $(0,2)$ tensor field on $M$. We define the {\it deviated $f$-energy
density} $E_{f,g,T}(\va)$ of the map $\va$ relative to $g$ and $T$ by the formula
$$
E_{f,g,T}(\va (t)) = {1 \ov 2} h^{\al\beta} (t) f_{ij} (x(t)) x^i_\al (t) x^j_\beta (t)\leqno(3)
$$
$$
+ {1 \ov 2}h^{\al\beta} (t) g_{ij} (x(t)) [x^i_\al (t) - T^i_\al (t)][x^j_\be (t) - T^j_\be (t)].
$$

{\bf Definition 2.} A map $\va$ is called {\it $f$-ultra-potential map} if it is a critical point of the
energy functional $E_{f, g, T}$, i.e., an extremal of the Lagrangian
$$
L_2 = E_{f, g, T}(\va)(t) \sqrt{|h|},
$$
for all compactly supported variations. The map $\va$ is called {\it  generalized
ultra-potential map relative to $g$ and $T$} if there exists a $(0,2)$-tensor field
$f$ on $(M,g)$  such that $\va$ is $f$-ultra-potential.

The {\it $f$-ultra-potential map equation} is a system of nonlinear {\it ultra-hyperbolic-Poisson PDEs}
and is expressed locally by
$$
h^{\al\be} x_{i\al\be} = h^{\al\be}g_{ij}\left[\frac{\pa T^j_{\al}}{\pa t^\be} - T^{j}_{\gamma}H^{\gamma}_{\al\be}\right]\leqno(4)$$
$$+h^{\al\be}x^j_{\al}T^k_{\be}\left[g_{si}G^{s}_{jk}-g_{sj}G^{s}_{ki}\right]+\frac{1}{2}T^j_{\al}T^k_{\be}\frac{\pa g_{jk}}{\pa x^i}.$$

Finally, if $g$ is a fixed Riemannian structure on $M$, let $X^i_\al (t,x)$ be a given
$C^\ty$ distinguished tensor field on $N \times M$ and
$c(t,x)$ be a given $C^\ty$ real function on $N\times M$. The {\it general energy
density} $E_{g,X}(\va)$ of the map $\va$, relative to $g$, $c$ and $X$ is defined by
$$
E_{g,c,X}(\va (t)) = {1 \ov 2} h^{\al\beta} (t) g_{ij} (x(t)) x^i_\al (t) x^j_\beta (t)\leqno(5)
$$
$$
- h^{\al\beta} (t) g_{ij} (x(t)) x^i_\al (t) X^j_\beta (t, x(t)) + c(t,x).
$$
Of course $E_{g,c,X}(\va)$ is a perfect square and is denoted by $E_{g,X}(\va)$ iff
$$
c = {1 \ov 2} h^{\al\beta}(t) g_{ij} (x(t)) X^i_\al (t, x(t))
X^j_\beta (t, x(t)).
$$
Similarly, for a relatively compact domain $\Omega \subset N$, we define the energy
$$
E_{g,X}(\va; \Omega) = \int_\Omega E_{g,X}(\va) (t) dv_h.
$$

{\bf Definition 3.} A map $\va$ is called {\it potential map} if it is a critical point of the
energy functional $E_{g,X}$, i.e., an extremal of the Lagrangian
$$
L_3 = E_{g,X}(\va)(t) \sqrt{|h|},
$$
for all compactly supported variations.

The {\it potential map equation} is a system
of nonlinear {\it elliptic-Poisson PDEs}, locally expressed  by
$$
 h^{\al\be} x^i_{\al\be} = g^{ij} {\pa c \ov \pa x^j} +
h^{\al\be} (\na_k X^i_\be
 - g_{kj} g^{il} \na_l X^j_\be)x^k_\al +
h^{\al\be} D_\al X^i_\be, \leqno (6)
$$
where $D$ is the covariant derivative on $(N,h)$ and $\na$ is the
covariant derivative on $(M,g)$. Explicitly, we have
$$
\na_j X^i_\al = {\pa X^i_\al \ov \pa x^j} + G^i_{jk} X^k_\al, \quad
D_\beta X^i_\al = {\pa X^i_\al \ov \pa t^\beta} - H^\gamma_{\beta\al} X^i_\gamma,
\leqno (7)
$$
$$
F_j{}^i{}_\al = \na_j X^i_\al - g_{hj} g^{ik} \na_k X^h_\al, \leqno (8)
$$
$$
{\pa g_{ij} \ov \pa x^k} = G^h_{ki} g_{hj} + G^h_{kj} g_{hi},
{\pa h^{\al\beta} \ov \pa t^\gamma} = - H^\al_{\gamma\lambda} h^{\la \beta}
- H^\beta_{\gamma\la} h^{\al\la}. \leqno (9)
$$

\section{Geometric dynamics and potential maps}

Let (N,h) and $(M,g)$ be two Riemannian manifolds of dimensions $m$,
respectively $n$ and let $X=(X^i_\al(t,x))$ be a $C^\ty$
distinguished tensor field on the manifold $N\times M$. The
classical geometric dynamics [7]-[13], [14]-[16] consists in
extending the normal PDE systems of first order
$$\frac{\pa x^j}{\pa t^\al}(t)=X^i_{\al}(t,x(t))\leqno(10)$$
into second order Euler-Lagrange type systems such that the solutions of the
system $(10)$ to be potential or harmonic maps relative to a certain geometric
structure. Following this idea, we recall, without proof, one of the main results in [7].

{\bf Theorem 1.} {\it Each solution $x: N\rightarrow M$ of the nonlinear and
non-homogeneous PDE system $(10)$ is a potential map. More precisely,
$x(\cdot)$ is an extremal for the least square type Lagrangian
$$L_4 =\di{1 \ov 2} h^{\al\beta} g_{ij} (x^i_\al - X^i_\al)(x^j_\beta - X^j_\beta)
\sqrt{|h|}.\leqno(11)$$}

\subsection {Geometric dynamics induced by \\non-homogeneous first order PDEs}

We start with a Riemannian manifold $(N,h)$ of dimension $m$, a
differential manifold $M$ of dimension $n$, a $C^{\ty}$ tensor field $Y=(Y^i_j(x))$ on
$M$, respectively a $C^{\ty}$ distinguished tensor field  $T=(T^i_\al(t))$ on $N$
and the implicit non-homogeneous nonlinear PDE system of order one

$$\frac{\pa x^j}{\pa t^\al}(t)Y^i_j(x(t))=T^i_{\al}(t).\leqno(12)$$

The purpose of this sub-section is to analyze the dynamics induced by the PDE
system $(12)$ and by appropriate metric tensor fields on $N$ and $M$. By differentiating the foregoing relation on $N\times M$ along a solution we obtain
$$\frac{\pa^2 x^j}{\pa t^\al \pa t^\be}Y^i_{j} + x^j_{\al}x^k_{\be}\frac{\pa Y^i_{j}}{\pa x^k} = \frac{\pa T^i_{\al}}{\pa t^\be}.$$
Using (12), adding-subtracting convenient terms, we change this relation into
$$\frac{\pa^2 x^j}{\pa t^\al \pa t^\be}Y^i_{j}-x^j_{\gamma}Y^i_{j}H^\gamma_{\al\be}+x^j_{\al}x^k_{\be}\frac{\pa Y^i_{j}}{\pa x^k}+x^j_{\al}x^k_{\be}G^i_{ks}Y^s_j=
\frac{\pa T^i_{\al}}{\pa t^\be}+T^j_{\al}x^{k}_{\be}G^{i}_{jk}-T^{i}_{\gamma}H^{\gamma}_{\al\be}.\leqno(13)$$
Taking the trace in $(13)$ with respect to $h^{\al\be}$, followed by a contraction with $g_{ip}$ and
adding the terms $-h^{\al\be}g_{sj}x^j_{\al}x^k_{\be}Y^i_{k}G^s_{ip}$ and $\frac{1}{2}x^j_{\al}x^k_{\be}Y^s_{j}Y^i_k\frac{\pa g_{is}}{\pa x^p}$,
we find

$$h^{\al\be}g_{ip}Y^i_{j}\left[\frac{\pa^2 x^j}{\pa t^\al \pa t^\be}-x^j_{\gamma}H^\gamma_{\al\be}\right]+h^{\al\be}x^j_{\al}x^k_{\be}\left[g_{ip}\frac{\pa Y^i_{j}}{\pa x^k}+g_{ip}G^i_{ks}Y^s_j\right.\leqno(14)$$
$$\left.-g_{sj}Y^i_{k}G^s_{ip}+\frac{1}{2}Y^j_{s}Y^k_i\frac{\pa g_{is}}{\pa x^p}\right]=h^{\al\be}g_{ip}\left[\frac{\pa T^i_{\al}}{\pa t^\be}-T^{i}_{\gamma}H^{\gamma}_{\al\be}\right]$$
$$+h^{\al\be}x^j_{\al}T^k_{\be}\left[g_{sp}G^{s}_{jk}-g_{sj}G^{s}_{kp}\right]+\frac{1}{2}T^j_{\al}T^k_{\be}\frac{\pa g_{jk}}{\pa x^p}.$$

If $$\Omega_{jk|p}=g_{ip}\frac{\pa Y^i_{j}}{\pa x^k}+g_{ip}G^i_{ks}Y^s_j-g_{sj}Y^i_{k}G^s_{ip}+\frac{1}{2}Y_j^{s}Y_k^i\left(g_{ij}G^j_{sp}+g_{sj}G^j_{ip}\right)$$
and $$S_{jk|p}=\frac{1}{2}\left[\Omega_{jk|p}+\Omega_{kj|p}\right],\,\,A_{jk|p}=\frac{1}{2}\left[\Omega_{jk|p}-\Omega_{kj|p}\right],$$
then
$$S_{jk|p}=\frac{1}{2}\left[g_{sp}\frac{\pa Y^s_j}{\pa x^k}+g_{sp}\frac{\pa Y^s_k}{\pa x^j}+Y^{s}_{j}\left(g_{ip}G^i_{ks}-g_{ik}G^i_{sp}\right)\right.$$
$$\left.+Y^{s}_{k}\left(g_{ip}G^i_{js}-g_{ij}G^i_{sp}\right)+Y_j^{s}Y_k^i\left(g_{ij}G^j_{sp}+g_{sj}G^j_{ip}\right)\right],$$

$$=\frac{1}{2}\left[g_{sp}\frac{\pa Y^s_j}{\pa x^k}+g_{sp}\frac{\pa Y^s_k}{\pa x^j}+Y^{s}_{j}\left(\frac{\pa g_{ps}}{\pa x^k}-\frac{\pa g_{ks}}{\pa x^p}\right)\right.$$
$$\left.+Y^{s}_{k}\left(\frac{\pa g_{ps}}{\pa x^j}-\frac{\pa g_{js}}{\pa x^p}\right)+Y_j^{s}Y_k^i\frac{\pa g_{is}}{\pa x^p}\right],$$
$$A_{jk|p}=\frac{1}{2}\left\{g_{sp}\left(\frac{\pa Y^s_j}{\pa x^k}-\frac{\pa Y^s_k}{\pa x^j}\right)+Y^s_{j}\frac{\pa g_{pk}}{\pa x^s}-Y^s_{k}\frac{\pa g_{pj}}{\pa x^s}\right\}.$$

If the $m$-sheet $x(\cdot)$ satisfies the PDE system $(12)$, then, along $x(\cdot)$
we have $$h^{\al\be}x^j_{\al}x^k_{\be}A_{jk|p}=0$$ and the PDE system  $(14)$ becomes
$$h^{\al\be}g_{ip}Y^i_{j}\left[\frac{\pa^2 x^j}{\pa t^\al \pa t^\be}-x^j_{\gamma}H^\gamma_{\al\be}\right]+h^{\al\be}x^j_{\al}x^k_{\be}S_{jk|p}
=h^{\al\be}g_{ip}\left[\frac{\pa T^i_{\al}}{\pa t^\be}-T^{i}_{\gamma}H^{\gamma}_{\al\be}\right]\leqno(15)$$
$$+h^{\al\be}x^j_{\al}T^k_{\be}\left[g_{sp}G^{s}_{jk}-g_{sj}G^{s}_{kp}\right]+\frac{1}{2}T^j_{\al}T^k_{\be}\frac{\pa g_{jk}}{\pa x^p}.$$

{\bf Theorem 2.} {\it The solutions of the implicit PDE system of first order $(12)$
are $f$-potential maps on $(M,g)$, relative to the distinguished tensor field $T$,
where the tensor field $f$ is solution for the PDE system
$$
\di{Y^s_i\frac{\pa f_{sj}}{\pa x^k}=\frac{\pa Y^s_j}{\pa x^k}f_{si}},\leqno(16)$$
satisfying also the conditions
$$\left\{\begin{array}{ll}f_{ij}=g_{is}(Y^s_j-\delta^s_j)\\
g_{is}Y^s_j=g_{js}Y^s_i.\end{array}\right.\leqno(17)$$

More precisely, the solutions of the implicit PDE system of first order $(12)$ are extremals for the Lagrangian
$$L_7=\frac{1}{2}h^{\al\be}\left[f_{ij}x^i_\al x^j_\be+g_{ij}(x^i_{\al}-T^i_{\al})(x^j_{\be}-T^j_{\be})\right]\sqrt{|h|}.$$}

{\bf Proof.} We consider first the Lagrangian
$$L_5 = h^{\al\be} \left(g_{ij} x^i_\al T^j_\be-\frac{1}{2}g_{ij} T^i_\al T^j_\be\right)
\sqrt{|h|}.
$$

In general, if $L = E \sqrt{|h|}$, where $E$ denotes an energy
density, then the Euler-Lagrange equations of extremals,
$$
{\pa L\ov \pa x^k} - {\pa \ov \pa t^\al} {\pa L\ov \pa x^k_\al} = 0
$$
can be written in the form
$$
{\pa E \ov \pa x^k} - {\pa \ov \pa t^\al}{\pa E \ov \pa x^k_\al} -
H^\gamma_{\gamma \al} {\pa E \ov \pa x^k_\al} = 0. \leqno (18)
$$
We compute
$$
\di{\pa E_5 \ov \pa x^k} = h^{\al\be} \di{\pa g_{ij} \ov \pa x^k}
x^i_\al T^j_\be-\frac{1}{2}h^{\al\be} \di{\pa g_{ij} \ov \pa x^k}
T^i_\al T^j_\be; $$

\hspace{1 cm} $$\di{\pa E_5 \ov \pa x^k_\al} = h^{\al\be}g_{kj} T^j_\be;$$

$$
- \di{\pa \ov \pa t^\al}\di{\pa E_5 \ov \pa x^k_\al} =
 - \di{\pa h^{\al\be} \ov \pa t^\al} g_{kj} T^j_\be -h^{\al\be} \di{\pa g_{kj} \ov \pa x^i} x^i_\al T^j_\be -h^{\al\be} g_{kj} \di{\pa T^j_\be \ov\pa t^\al}.
$$
Replacing in $(18)$, we find
$$
-\delta L_5=h^{\al\be}g_{ik}\left[\frac{\pa T^i_{\be}}{\pa t^\al}-T^{i}_{\gamma}H^{\gamma}_{\al\be}\right]
+h^{\al\be}x^i_{\be}T^j_{\be}\left[g_{sk}G^{s}_{ij}-g_{si}G^{s}_{jk}\right]+\di{1 \ov 2} h^{\al\be} \di{\pa g_{ij} \ov \pa x^k} T^i_\al T^j_\be,
$$
which is precisely the  right hand in $(15)$.

Next, we shall compute the first variation for the Lagrangian
$$L_6 = \di{1 \ov 4} h^{\al\be} (g_{is}Y^s_j+g_{js}Y^s_i) x^i_\al x^j_\be\sqrt{|h|}.$$

We obtain
$$
\di{\pa E_6 \ov \pa x^k} = \di{1 \ov 4} h^{\al\beta} \left[\di{\pa g_{is} \ov
\pa x^k} Y^s_j + g_{is} \di{\pa Y^s_j \ov \pa x^k} +  \di{\pa g_{js} \ov
\pa x^k} Y^s_i + g_{js} \di{\pa Y^s_i \ov \pa x^k}\right] x^i_\al x^j_\be,
$$
$$\di{\pa E_6 \ov \pa x^k_\al} = \di{1 \ov 2}h^{\al\beta}(g_{ks}Y^s_j+g_{js}Y^s_k) x^j_\beta, $$

$$
- \di{\pa \ov \pa t^\al}\di{\pa E_6 \ov \pa x^k_\al} = - \di{\pa h^{\al\beta}
\ov \pa t^\al} g_{ks}Y^s_j x^j_\beta -h^{\al\beta}
g_{ks}Y^s_j \frac{\pa^2 x^j}{\pa t^\al \pa t^\be}
$$
$$- \di{1 \ov 2}h^{\al\beta} \left[\di{\pa g_{ks} \ov
\pa x^i} Y^s_j + g_{ks} \di{\pa Y^s_j \ov \pa x^i} +  \di{\pa g_{js} \ov
\pa x^i} Y^s_k + g_{js} \di{\pa Y^s_k \ov \pa x^i}\right] x^i_\al x^j_\be.$$

Replacing in $(18)$, we find
$$-\delta L_6=h^{\al\be}g_{sk}Y^s_{j}\left[\frac{\pa^2 x^j}{\pa t^\al \pa t^\be}-
x^j_{\gamma}H^\gamma_{\al\be}\right]-h^{\al\be}x^i_{\al}x^j_{\be}\Sigma_{ij|k},$$
where
$$\Sigma_{ij|k}=\frac{1}{4}\left[Y^s_j\left(\frac{\pa g_{is}}{\pa x^k}-\frac{\pa g_{ks}}{\pa x^i}\right)+
Y^{s}_{i}\left(\frac{\pa g_{js}}{\pa x^k}-\frac{\pa g_{ks}}{\pa x^j}\right)-g_{sk}\left(\frac{\pa Y^s_i}{\pa x^j}+\frac{\pa Y^s_j}{\pa x^i}\right)\right.
$$
$$\left. -Y^{s}_{k}\left(\frac{\pa g_{js}}{\pa x^i}+\frac{\pa g_{is}}{\pa x^j}\right)+g_{sj}\left(\frac{\pa Y^s_i}{\pa x^k}-
\frac{\pa Y^s_k}{\pa x^i}\right)+g_{si}\left(\frac{\pa Y^s_j}{\pa x^k}-\frac{\pa Y^s_k}{\pa x^j}\right)\right].$$
By computation, using relations $(17)$, we obtain
$$\Sigma_{ij|k}+S_{ij|k}=\frac{1}{2}\left(Y^s_i\frac{\pa f_{sj}}{\pa x^k}-\frac{\pa Y^s_j}{\pa x^k}f_{si}\right)$$
that is, using relation $(16)$,
$$\Sigma_{ij|k}=-S_{ij|k}.$$
Therefore, $-\delta L_6$ has the same expression as the left hand side in relation $(15)$.
We conclude that $x(\cdot)$ is an extremal for the Lagrangian
$$L_7  =L_6 - L_5$$
$$ =\di{1 \ov 4} h^{\al\be} (g_{is}Y^s_j+g_{js}Y^s_i) x^i_\al x^j_\be\sqrt{|h|}-
h^{\al\be} \left(g_{ij} x^i_\al T^j_\be-\frac{1}{2}g_{ij} T^i_\al T^j_\be\right)\sqrt{|h|}$$
$$
= \frac{1}{2}h^{\al\be}\left[f_{ij}x^i_\al x^j_\be+g_{ij}(x^i_{\al}-T^i_{\al})(x^j_{\be}-T^j_{\be})\right]\sqrt{|h|}.
$$\qed

\subsection {Geometric dynamics induced by \\homogeneous first order PDEs}

Now, let us consider the homogeneous nonlinear first order PDE system
$$\frac{\pa x^j}{\pa t^\al}(t)Y^i_j(x(t))=0.\leqno(19)$$
By differentiating the tensor field $x^j_{\al} Y^i_j(x(t))dt^\al\otimes\frac{\pa}{\pa x^i}$ on
$N\times M$ along a solution and adding-subtracting appropriate terms, we obtain
$$\frac{\pa^2 x^j}{\pa t^\al \pa t^\be}Y^i_{j}-x^j_{\gamma}Y^i_{j}H^\gamma_{\al\be}+x^j_{\al}x^k_{\be}\nabla_kY^i_j+x^j_{\al}x^k_{\be}G^s_{kj}Y^i_s=0,$$
or
$$Y^i_{j}\left[\frac{\pa^2 x^j}{\pa t^\al \pa t^\be}-x^j_{\gamma}H^\gamma_{\al\be}+x^l_{\al}x^k_{\be}G^j_{kl}\right]+x^j_{\al}x^k_{\be}(\nabla_kY)^i_j=0,\leqno(20)$$
where $$(\nabla_kY)^i_j=\frac{\pa Y^i_{j}}{\pa x^k}+G^i_{ks}Y^s_j-G^s_{kj}Y^i_s.$$

Taking the trace of $(20)$ with respect to $h^{\al\be}$ and lowering the index $i$
with $g_{ik}$, we get
$$h^{\al\be}g_{ik}Y^i_{j} x^j_{\al\be}+h^{\al\be}g_{ik}x^j_{\al}x^p_{\be}\nabla_pY^i_j=0,\leqno(21)$$
where $$x^i_{\al\be}=\frac{\pa^2 x^i}{\pa t^\al \pa t^\be}-x^i_{\gamma}H^\gamma_{\al\be}+x^j_{\al}x^{k}_{\be}G^i_{jk}.$$

{\bf Theorem 3.} {\it The solutions of the implicit homogeneous PDE system of first
order $(19)$ are $f$-harmonic maps on $M$, where $f_{ij}=g_{is}Y^s_{j}$ and $g$ is  solution for the PDE system
$$(\nabla_{k}Y)^i_j=\frac{\pa Y^i_{j}}{\pa x^k}+G^i_{ks}Y^s_j-G^s_{kj}Y^i_s=0\leqno(22)$$
satisfying the symmetry condition
$$g_{is}Y^s_j=g_{js}Y^s_i.\leqno(23)$$
Here $G^k_{ij}$ mean the Christoffel symbols of $g$.

Moreover, the solutions of the implicit homogeneous PDE system of first order $(19)$ are extremals for the Lagrangian
$$L_8=\frac{1}{2}h^{\al\be}f_{ij}x^i_\al x^j_\be\sqrt{|h|}.$$}

{\bf Remarks.} (1) The idea of finding $G^i_{ks}$ from the relation
(22) was developed in [17].

(2) Writing the complete integrability conditions for the PDEs $(22)$, we obtain
$$Y^i_s R^s_{jkl}=Y^s_j R^i_{skl},\leqno(24)$$
where $R$ denotes the Riemann curvature tensor field corresponding to the solution $g$.

{\bf Proof.} We need to verify that the PDE system $(21)$ is in fact the Euler-Lagrange PDE system corresponding to the Lagrangian
$$L_8=\frac{1}{2}h^{\al\be}g_{ik}Y^k_jx^i_\al x^j_\be\sqrt{|h|}.$$
On the other hand, we know that
$$-\delta L_8=h^{\al\be}g_{sk}Y^s_{j}\left[\frac{\pa^2 x^j}{\pa t^\al \pa t^\be}-
x^j_{\gamma}H^\gamma_{\al\be}\right]-h^{\al\be}x^i_{\al}x^j_{\be}F_{ij|k},$$
and the hypotheses ensure us that $F_{ij|k}=-g_{ks}Y^s_lG^l_{ij}$. We obtain
$$-\delta L_8=h^{\al\be}g_{sk}Y^s_{j}\left[\frac{\pa^2 x^j}{\pa t^\al \pa t^\be}-
x^j_{\gamma}H^\gamma_{\al\be}+x^i_{\al}x^l_{\be}\Gamma^j_{il}\right]=-h^{\al\be}g_{sk}Y^s_{j}x^j_{\al\be}$$
and the Euler Lagrange PDE system corresponding to $L_8$ has the same expression as in $(21)$.
\qed

\section{Gauss equations for an \\integral submanifold map}

In this section, $(N,h)$ and $(M,g)$ denote an $m$-dimensional, respectively,
an $n$-dimensional Riemannian manifold and $X = X^i_{\al}(t, x)dt^\al\otimes\frac{\pa}{\pa x^i}$
is a $C^\ty$ distinguished tensor field on $N\times M$, satisfying the integrability conditions
$$\frac{\pa X^i_{\al}}{\pa x^j}X^j_{\be}=\frac{\pa X^i_{\be}}{\pa x^j}X^j_{\al}.$$
We are looking for describing the geometry of the $C^\ty$ integral submanifolds
$$x: N\rightarrow M,\,\,\frac{\pa x^i}{\pa t^\al}(t)=X^i_{\al}(t,x(t)).\leqno(25)$$

Differentiating PDEs $(25)$ along a solution and replacing $x^j_{\be} = X^j_{\be}$, we find
$$\frac{\pa^2 x^i}{\pa t^\al \pa t^\be}=\frac{\pa X^i_{\al}}{\pa x^j}X^j_{\be}+\frac{\pa X^i_{\al}}{\pa t^\be}.\leqno(26)$$
On the other side, the Gauss equation corresponding to an $m$-dimensional submanifold $x(\cdot)$ is of the form
$$\frac{\pa^2 x^i}{\pa t^\al \pa t^\be}(t)=\Lambda^\gamma_{\al\be}(t)x^i_{\gamma}(t)+\Lambda^a_{\al\be}(t)N^i_{a}(x(t)),\leqno(27)$$
where $N_a=N_a^i\frac{\pa}{\pa x^i}$ is an orthonormal family of vector fields on $M$, normal to the submanifold $x(N)$, that is
$$g_{ij}N^i_a N^j_b=\delta_{ab},\,\,g_{ij}N^i_a X^j_\be=0.\leqno(28)$$
Moreover, let $h_{\al\be}(t)=\left(g_{ij}X^i_{\al}X^j_{\be}\right)(x(t))$.
From the relations (26)-(27), we obtain the Tzitzeica connection
$$\Lambda^\gamma_{\al\be}(t)=h^{\gamma\sigma}(t)g_{ik}X^k_{\sigma}\left[\frac{\pa X^i_{\al}}{\pa x^j}X^j_{\be}+
\frac{\pa X^i_{\al}}{\pa t^\be}\right](x(t)),\leqno(29)$$
and the fundamental forms
$$\Lambda^a_{\al\be}(t)=\delta^{ab}(t)g_{ik}N^k_{b}\left[\frac{\pa X^i_{\al}}{\pa x^j}X^j_{\be}+\frac{\pa X^i_{\al}}{\pa t^\be}\right](x(t)).\leqno(30)$$

\section{General potentiality of submanifold maps}

Our aim is to prove that there exists an infinity of Riemannian structures such that the lift of a
submanifold map to the jet bundle of order one is a potential map.
Let $x: N\rightarrow M$ be a $C^\ty$ $m$-dimensional Riemannian submanifold of $(M,g)$ . Then, the Gauss formula of the submanifold $x$ is
$$\frac{\pa^2 x^i}{\pa t^\al\pa t^\be}=\Lambda^{\gamma}_{\al\be}\, x^i_{\gamma} + \Lambda^{a}_{\al\be}N^i_a,\leqno(Gauss)$$
where $\{N_s|s = 1, ..., n-m\}$ denotes a family of normal vector fields to the submanifold, $\Lambda^{\gamma}_{\al\be}$
are the components of the connection and $\Lambda^{a}_{\al\be}$ are the fundamental forms.
We make the assumption that $\{N_s|s = 1, ..., n-m\}$ is an orthonormal distribution.
We transform the Gauss second order PDE system into a first order system on the jet bundle $J^1(N,M)$ as it follows:
$$ \left\{\begin{array}{ll}
\frac{\pa x^i}{\pa t^{\gamma}}&=x^i_\gamma,\\
\frac{\pa x^i_\al}{\pa  t^\be}&=\Lambda^{\gamma}_{\al\be}x^i_{\gamma}+\Lambda^a_{\al\be}N^i_a.\end{array}\right.  \leqno(Gauss)$$

Let $\eta$ be the induced Riemannian metric on the submanifold $N$, i.e.,
$$\eta_{\alpha\beta}(t)=g_{ij}(x(t))x^i_\al(t)x^j_\be(t).$$
Moreover,
$$\Lambda^{\gamma}_{\al\be}=\frac{1}{2}\eta^{\gamma\sigma}\left[\frac{\pa \eta_{\al\sigma}}{\pa t^\beta}+
\frac{\pa \eta_{\be\sigma}}{\pa t^\al}-\frac{\pa \eta_{\al\be}}{\pa t^\sigma}\right];\,\,\Lambda^{a}_{\al\be}=
g_{ij}\frac{\pa^2 x^i}{\pa t^\al \pa t^\be}N^j_b\delta^{ab}$$
denote the Christoffel symbols, respectively the second fundamental forms of the submanifold.

Let $h$ be an arbitrary Riemannian structure on $N$ and let $J^1(N,M)$, with local coordinates
$(t^\al,\x^i_0=x^i,\x^i_{\al}=x^i_{\al})$, denote the first order jet bundle.
Let ${\mathcal I}=\{I=\left(^i_\al\right)|i=1,..,n,\,\,\al=0,...,m\}$ and
${\bf \varphi}:N\rightarrow J^1(N,M),$ ${\bf\varphi}(t)=(t,\x^I(t))$. Then, we may write the Gauss second order PDE system
as the Gauss first order PDE system in the jet bundle of order one,
 $$
\frac{\pa \x^I}{\pa t^{\mu}}(t)=X^I_\mu(t,\x(t)), \leqno(Gauss)$$
where
$$X^I_\mu(t,\x)=\left\{\begin{array}{ccc}
x^i_\mu, & if & I=\left(^i_0\right)\\
\Lambda^{\gamma}_{\mu\al}x^i_{\gamma}+\Lambda^{a}_{\mu\al}N^i_a,& if &I=\left(^i_\al\right), \,\,\al\neq0.
\end{array}\right.\leqno(31)$$

We know from [7] that the solutions of a normal system of PDEs of order one are potential
maps in an appropriate geometrical structure. The purpose of this paper is to prove that,
for each embedded submanifold, there are geometric structures on the environmental
manifold such that the lift to the jet bundle of a submanifold map is a potential map
and to find the PDEs describing this Riemannian structures.
Let $h_{\al\be}(t)dt^\al\otimes dt^\be+\g_{IJ}d\x^I\otimes d\x^J$ be an arbitrary Riemannian
structure on $J^1(N,M)$. The following result is a consequence of Theorem 1.

{\bf Theorem 4.} {\it The lift of a submanifold map $x:N\rightarrow M$ to the jet bundle $J^1(N,M)$
is a potential map. More precisely, it is an extremal for all the least squares Lagrangians (depending on the Riemannian structure $\g$)
$$L_{\g} = \di{1 \ov 2} h^{\mu\nu} \g_{IJ} (\x^I_\mu - X^I_\mu)(\x^J_\nu - X^J_\nu)
\sqrt{|h|} $$
$$=\di{1 \ov 2} h^{\mu\nu} \g_{\left(^i_{\al}\right)\left(^j_{\be}\right)} \left(\x^i_{\mu\al} -
 \Lambda^{\gamma}_{\mu\al}x^i_{\gamma}-\Lambda^{a}_{\mu\al}N^i_a\right)\left(\x^j_{\nu\be} -
\Lambda^{\gamma}_{\nu\be}x^j_{\gamma}-\Lambda^{a}_{\nu\be}N^j_a\right)
\sqrt{|h|}. $$}

{\bf Remarks.} (1) Writing the Euler-Lagrange PDEs for the Lagrangian $L_{\g}$, we obtain
$$
h^{\mu\nu} \x^I_{\mu\nu} = \g^{IL} h^{\mu\nu} \g_{KJ} (\na_L X^K_\mu)
X^J_\nu + h^{\mu\nu} F_J{}^I{}_\mu \x^J_\nu + h^{\mu\nu} D_\nu X^I_\mu,
\leqno(E-L)_\g
$$
where
$$\na_L X^K_\mu=\frac{\pa X^K_\mu} {\pa \x^L} +\Gamma^K_{LS}X^S_\mu,\,\, D_\nu X^I_\mu = - H^\gamma_{\mu\nu} X^I_\gamma,$$
$$F_J{}^I{}_\mu = \na_J X^I_\mu - \g^{IL} \g_{KJ} \na_L X^K_\mu,$$
and
$$\x^I_{\mu\nu} =
{\pa^2 \x^I \ov \pa t^\mu \pa t^\nu} -
H^\gamma_{\mu\nu} \x^I_\gamma + \Gamma^I_{JK} \x^J_\mu \x^K_\nu.
$$
(2) There exists an infinity of geometrical structures $\g$ such
that the lift of a submanifold map is a potential map.

\section{General harmonicity of submanifold maps}

Generally, an arbitrary submanifold map between two Riemannian
manifolds $(N,g)$ and $(M,g)$ is not a harmonic one and not even a
potential one. Nevertheless, Theorem 4 proved that its lift to the
first order jet bundle, endowed with an infinite possible Riemannian
structures, it is a potential map. We shall see further, that the
submanifold map may also be harmonic, in a general sense. Indeed,
let $N$ and $M$ be two differentiable manifolds and $x:N\rightarrow
M$ be a differentiable submanifold map. Let $\nabla x(t)=(\frac{\pa x^i}{\pa t^\al}(t))$
be the Jacobian matrix  which is of rank $m$. For each point $t\in N$, the algebraic system
$$\frac{\pa x^j}{\pa t^\al}(t)\xi^i_j(t)=0,$$
defines the matrix function $\xi^i_j(t)$. Let $Y=Y^i_j(x)dx^j\otimes \frac{\pa}{\pa x^i}$ be a
$C^{\ty}$ tensor field on $M$ such that $Y^i_j(x(t))=\xi^i_j(t).$
Then $x(\cdot)$ is a solution for the nonlinear homogeneous PDE system
$$\frac{\pa x^j}{\pa t^\al}(t)Y^i_j(x(t))=0.\leqno(32)$$
As a consequence of Theorem 3, an arbitrary Riemannian structure $h$ on $N$,
together with a Riemannian structure $g$ on $M$ (solution for a nonlinear PDE
system marking a parallelism condition, satisfying also a symmetry condition)
determine the general harmonicity of the map $x(\cdot)$.

In the sequel, we shall describe an alternative way of obtaining general harmonicity,
where the Riemannian structure $h$ stays fixed, but the symmetry condition for the
Riemannian structure $g$ is unconditional.

For this, we start from the relation
$$h^{\al\be}g_{is}Y^s_{j}\left[\frac{\pa^2 x^j}{\pa t^\al \pa t^\be}-
x^j_{\gamma}H^\gamma_{\al\be}\right]+h^{\al\be}g_{is}x^p_{\al}x^k_{\be}\left[\frac{\pa Y^s_p}{\pa x^k}+Y^j_pG^s_{jk}\right]=0,\leqno(33)$$
obtained by differentiating the initial homogeneous system $(32)$
along a solution and taking, afterwards, the trace with respect to
$h$ and the contraction with respect to $g$. We know that, for a
fixed Riemannian structure $g^0$ on $M$ and a family of normal
vector fields $N_a=(N^i_a),\,\,a=1,...,n-m$, each submanifold map $x: N\rightarrow M$ satisfies the Gauss equations
$$\frac{\pa^2 x^i}{\pa t^\al\pa t^\be}=\Lambda^{\gamma}_{\al\be}x^i_\gamma+\Lambda^{a}_{\al\be}N^i_a,\leqno(Gauss)$$
where, if $h$ is the metric induced  by $g^0$ on $N$, then $\Lambda^{\gamma}_{\al\be}=H^{\gamma}_{\al\be}$
are the components of the corresponding Levi-Civita connection.
Let us choose  $\Lambda^{0\gamma}_{\al\be}=H^{0\gamma}_{\al\be}$ such that
$$\frac{\pa^2 x^i}{\pa t^\al\pa t^\be}=\Lambda^{0\gamma}_{\al\be}x^i_\gamma$$
and let $h^0$ be solution for the Ricci PDE system
$$\frac{\pa h^0_{\al\be}}{\pa t^\gamma}=h^0_{\al\sigma}\Lambda^{0\sigma}_{\be\gamma}+h^0_{\be\sigma}\Lambda^{0\sigma}_{\al\gamma}.$$

Using this particular structure allows us to replace relation $(33)$ with
$$h^{0\al\be}\left(g_{is}Y^s_{j}+g_{js}Y^s_{i}\right)\left[\frac{\pa^2 x^j}{\pa t^\al \pa t^\be}-
x^j_{\gamma}H^{0\gamma}_{\al\be}\right]+h^{0\al\be}g_{is}x^p_{\al}x^k_{\be}\left[\frac{\pa Y^s_p}{\pa x^k}+
Y^j_pG^s_{jk}\right]=0.\leqno(34)$$
Let $\Omega_{pk|i}=g_{is}\left[\frac{\pa Y^s_p}{\pa x^k}+Y^j_pG^s_{jk}\right]$ and $S_{pk|i}=
\frac{1}{2}(\Omega_{pk|i}+\Omega_{kp|i})$, $A_{pk|i}=\frac{1}{2}(\Omega_{pk|i}-\Omega_{kp|i})$.
By computation, we obtain that each solution $x(\cdot)$ of PDE system $(32)$ satisfies the
equality $h^{0\al\be}x^p_{\al}x^k_{\be}A_{pk|i}=0$, and therefore, the relation $(34)$ becomes
$$h^{0\al\be}\left(g_{is}Y^s_{j}+g_{js}Y^s_{i}\right)\left[\frac{\pa^2 x^j}{\pa t^\al \pa t^\be}-
x^j_{\gamma}H^{0\gamma}_{\al\be}\right]+h^{0\al\be}x^p_{\al}x^k_{\be}S_{pk|i}=0.\leqno(35)$$

{\bf Theorem 5.} {\it The solutions of the implicit homogeneous PDE system of first order $(32)$ are
$f$-harmonic maps relative to $(N,h^0)$ and $M$, where $f_{ij}=g_{is}Y^s_{j}+g_{js}Y^s_{i}$ and $g$ is  solution for the PDE system
$$g_{is}\left[(\nabla_{k}Y)^s_j-(\nabla_{j}Y)^s_k\right]+g_{js}\left[(\nabla_{k}Y)^s_i-(\nabla_{i}Y)^s_k\right]=2g_{sp}G^p_{ij}Y^s_k.\leqno(36)$$}

{\bf Hint.} We consider the Lagrangian $L_9=\frac{1}{2}h^{0\al\be}f_{ij}x^i_\al x^j_{\be}\sqrt{h^0}$.
Similar arguments with the forgoing one ensure us that relations $(35)$ describe the Euler-Lagrange PDE system corresponding to this Lagrangian. \qed

Constantin Udriste\\
{\it University Politehnica of Bucharest\\ Faculty of Applied
Sciences\\ Department Mathematics-Informatics I\\ Splaiul
Independentei 313, 060042, Bucharest,
Romania,\\email: anet.udri@yahoo.com}\\

Vasile Arsinte\\ {\it Callatis High School, Rozelor 36, Mangalia,
Romania,\\ email: varsinte@seanet.ro}\\

Andreea Bejenaru \\
{\it University Politehnica of Bucharest\\ Faculty of Applied
Sciences\\ Department Mathematics-Informatics I\\ Splaiul
Independentei 313, 060042, Bucharest,
Romania,\\
 email: bejenaru.andreea@yahoo.com}

\end{document}